\documentclass[11pt]{article}

\title{\Large ON A STOCHASTIC WAVE EQUATION\\ DRIVEN BY A NON-GAUSSIAN L\'{E}VY PROCESS
\footnote{Supported by the LPMC at Nankai University and the NSF of
China (No.10871103, No.60874085). Email address:
bolijunnk@yahoo.com.cn (L. Bo), kehuashink@gmail.com (K. Shi, the
corresponding author), yjwang@nankai.edu.cn (Y. Wang)}}

\author{\small LIJUN BO$^{a}$, KEHUA SHI$^{b}$, YONGJIN WANG$^{b}$\\
\small\it{ $^{a}$Department of Mathematics, Xidian University, Xi'an
$710071$, P.R. China} \\
\small\it{ $^{b}$School of Mathematical Sciences, Nankai University,
Tianjin $300071$, P.R. China}}

{\footnotesize\date{\today}}

\usepackage{amsfonts}

\newcommand{\E}{\mathbf{E}}
\newcommand{\EE}{{\mathrm{E}}}

\newcommand{\N}{{\mathbf{N}}}
\newcommand{\R}{{\mathbf{R}}}
\newcommand{\Y}{\mathbf{P}}

\newcommand{\F}{{\mathcal{F}}}

\newcommand{\D}{{\mathrm{d}}}
\newcommand{\DD}{{\mathbb{D}}}
\newcommand{\DDD}{{{D}}}
\newcommand{\A}{{\mathbf{A}}}
\newcommand{\T}{{\mathrm{T}}}
\newcommand{\B}{{\mathcal{B}}}

\newcommand{\LL}{{{L}}}
\newcommand{\p}{{\mathcal{P}}}
\newcommand{\q}{{\mathrm{p}}}

\newcommand{\J}{{\mbox{-}}}
\newcommand{\prf}{\noindent{\bf Proof}.\quad}
\newcommand{\prfe}{{\hfill$\Box$}}

\newtheorem{remark}{Remark}[section]
\newtheorem{lemma}{Lemma}[section]

\newtheorem{theorem}{Theorem}[section]
\newtheorem{proposition}{Proposition}[section]

\begin{document}

\maketitle

\begin{abstract}
This paper investigates a damped stochastic wave equation driven by
a non-Gaussian L\'{e}vy noise. The weak solution is proved to exist
and be unique. Moreover we show the existence of a unique invariant
measure associated with the transition semigroup under mild conditions.\\

\noindent{Key words:}\quad   Damped wave equation, L\'{e}vy noise,
invariant measure\\
\noindent{MSC:}\quad 60H15; 35K90; 47D07

\end{abstract}

\section{Introduction}
\setcounter{equation}{0}

Let $(\Omega,\F,(\bar{\F}_t)_{t\geq0},\Y)$ be a complete filtered
probability space, and on which, $\widetilde{N}(\D z,\D t):=N(\D
z,\D t)-\pi(\D z)\D t$ defines a compensated Poisson random measure
of a Poisson random measure
$N:\B(Z\times[0,\infty))\times\Omega\to\N\cup\{0\}$ with the
characteristic measure $\pi(\cdot)$ on $(Z,\B(Z))$ with $Z=\R^m$
($m\in\N$). The characteristic measure $\pi(\cdot)$ satisfies that
\begin{eqnarray}
\pi(\{0\})=0,\ \ \ \ \ \ \ \ \int_{Z}1\wedge|z|^2\pi(\D z)<\infty.
\end{eqnarray}
According to (1.1), for $Z_1=\{z\in Z;\ |z|\leq1\}$, we can define
\begin{eqnarray}
\bar{\theta}=\int_{Z_1}|z|^2\pi(\D z),\ \ \ \ \ \ \ \ \
\underline{\theta}=\pi(Z\setminus Z_1).
\end{eqnarray}
 In the current paper, we
are concerned with the following hyperbolic equation with a
non-Gaussian L\'{e}vy noise perturbation:
\begin{eqnarray}
\left\{\begin{array}{ll} \frac{\partial^2u(t,\xi)}{\partial
t^2}+\kappa \frac{\partial u(t,\xi)}{\partial t}-\Delta u(t,\xi)\\ \
\ \ \ \ =\int_{Z_1}a(u(t\J,\xi),z)\dot{\widetilde{N}}(\D z,t)\\ \ \
\ \ \ \ \ \ \ \ +\int_{Z\setminus Z_1}b(u(t\J,\xi),z)\dot{N}(\D
z,t),
\ \ \ \ \ (t,\xi)\in[0,\infty)\times D,\\
u(0,\xi)=\varphi(\xi),\ \ \frac{\partial u(0,\xi)}{\partial t}=\psi(\xi),\ \ \ \xi\in D,\\
u(t,\xi)=0,\ \ \ \ \ \ (t,\xi)\in[0,\infty)\times\partial D,
\end{array}\right.
\end{eqnarray}
where the domain $D\subset \R^d$ is a bounded open set with
sufficiently regular boundary $\partial D$ and $\kappa>0$ denotes
the damped coefficient. The random measure $\widetilde{N}(\D z,\D
t)=N(\D z,\D t)-\pi(\D z)\D t$ denotes the compensated Poisson
random measure through the compensator of $N(\D z,\D t)$. In
addition, the functions $a:\R\times Z_1\to\R$ and $b:\R\times
Z\setminus Z_1\to\R$ are some regular functions with the exact
conditions in Section 2 below.

White noise perturbed stochastic wave equations have been
investigated in the literature (see e.g. \cite{BDT, BSW08, BTW, BMS,
Chow02, Chow06, PeZabc} and the references therein). In Chow
\cite{Chow02}, the global (weak) solutions of stochastic wave
equations with polynomial nonlinearity were explored by constructing
appropriate Lyapunov functionals. In a successive paper, Chow
\cite{Chow06} discussed the asymptotic behavior of the global (weak)
solution to a semilinear stochastic wave equation by using the
energy approach. Brzeniak et al. \cite{BMS} studied an abstract
stochastic wave equation: stochastic beam equation and Lyapunov
functions techniques were used to prove the existence of global mild
solutions and asymptotic stability of the zero solution. Barbu et
al. \cite{BDT} demonstrated the existence of an invariant measure
for the transition semigroup associated with a stochastic wave
equation with the nonlinear dissipative damping and further
established the uniqueness in some special case. In Bo et al.
\cite{BTW}, the authors used appropriate energy inequalities to give
sufficient conditions such that the local solutions of a class of
(strong) damped stochastic wave equations are blowup with a positive
probability or explosive in $L^2$-sense.

A recent work in Peszat and Zabczyk \cite{PeZabc06, PeZabc}  was to
consider the following wave equation driven by an impulsive noise,
\begin{eqnarray}
\frac{\partial^2 u(t)}{\partial t^2}=[\Delta u(t)+f(u(t))]\D
t+b(u(t))P\D Z(t),
\end{eqnarray}
where $f, b: \R\to\R$ are Lipschitz continuous, $P$ is a
regularizing linear operator and the impulsive noise
$Z=(Z_t)_{t\geq0}$ is formulated as a Poisson random measure. By
estimating the stochastic convolution w.r.t. Poisson random measure,
the authors proved that (1.4) admits a unique mild solution,
provided the intensity measure of $Z$ and eigenvectors of the
Laplace operator jointly satisfy a finite infinite series condition.

Compared with the above mentioned literature, we discuss several
other aspects of the differences in this article. First, the
objective equation we considered is the damped wave equation (with
the damped term $\kappa\frac{\partial u(t)}{\partial t}$) which is
used to model nonlinear phenomena in relativistic quantum mechanics
with local interaction (see e.g. \cite{Sch, Segal}). Second, this
paper focuses on the notion of the weak solution which is a stronger
form than the mild notion. Third, the perturbation can include a
general non-Gaussian L\'{e}vy noise which is much wider than the one
considered in \cite{PeZabc06, PeZabc}. Specially, we don't make any
assumptions for the L\'{e}vy measure in the process of proving the
existence and uniqueness of the weak solution. Finally, we also
explore the invariant measure associated with the weak solution,
which was not considered in \cite{PeZabc06, PeZabc}.

The paper is organized as follows: In the coming section, some
preliminaries and hypothesis are given. In Section 3, the existence
of a unique weak solution to (1.3) is established. Section 4 is
devoted to proving the existence of a unique invariant measure
corresponding to the weak solution under mild conditions.

\section{Preliminaries and hypothesis}
\setcounter{equation}{0}

We begin with some basic notation, functional spaces and
inequalities, which will be used frequently in the following
sections.

Define a linear operator $A$ by
\begin{eqnarray}
Au =-\Delta u,\ \ \ u\in\DDD(A)=H^2(D)\cap H^1_0(D).
\end{eqnarray}
where $H^p(D)$ is the set of all functions $u\in L^2(D)$ which have
generalized derivatives up to order $p$ such that $D^\alpha u\in
L^2(D)$ for all $\alpha:|\alpha|\leq p$, and $H^p_0(D)$ denotes the
closure of $C_0^\infty(D)$ in $H^p(D)$. Set $H=L^2(D)$ and
$V=H_0^1(D)$. Then $A$ is a positive self-adjoint unbounded operator
on $H$. On the other hand, both $H$ and $V$ are Hilbert spaces if we
endow them with usual inner products $\left<\cdot,\cdot\right>$ and
$\ll\cdot,\cdot\gg$, respectively. Furthermore,
\begin{eqnarray}
\DDD(A)\subset V\subset H\subset V^*,
\end{eqnarray}
where $V^*$ denotes the dual space of $V$, and the embedding
$V\subset H$ is compact. Thus there exists an orthonormal basis of
$H$, $(e_k)_{k=1,2,\dots}$ which consists of eigenvectors of $A$
such that $Ae_k=\lambda_ke_k$ for $k=1,2,\dots$ and
$0<\lambda_1\leq\lambda_2\leq\dots,$ with
$\lim_{k\to\infty}\lambda_k=+\infty$. According to the spectral
theory, for each $s\in\R$, we can define Hilbert space
$V_{2s}=\DDD(A^s)$, under the following inner product and the norm:
\begin{eqnarray}
\left<u,v\right>_{2s}&:=&\sum_{k=1}^\infty\lambda_k^{2s}\left<u,e_k\right>\left<v,e_k\right>,\\
|u|_{2s}&:=&\left[\sum_{k=1}^\infty\lambda_k^{2s}\left|\left<u,e_k\right>\right|^2\right]^{1/2}.
\end{eqnarray}
Obviously $\left<\cdot,\cdot\right>=\left<\cdot,\cdot\right>_0$ and
$\ll\cdot,\cdot\gg=\left<\cdot,\cdot\right>_1$. For parsimony, we
set $|\cdot|=|\cdot|_0$ and $\|\cdot\|=|\cdot|_1$. The following
Poincare-type inequality are well known (see e.g. Temam \cite{Temam}
and Zeidler \cite{Zei}):
\begin{eqnarray}
|u|_{\alpha_1}&\leq&\lambda_1^{\frac{\alpha_1-\beta_1}{2}}|u|_{\beta_1},\
\ {\rm for}\ \alpha_1\leq\beta_1,\ {\rm and}\ u\in
\DDD(A^{\beta_1/2}).
\end{eqnarray}

At the end of this section, we make the following basic assumptions:\\

\noindent$({H1})$\quad $a,b:\R\times Z\to \R$ are measurable and
there exists a constant $\ell_a>0$ such that
\begin{eqnarray*}
a(0,z)&\equiv& 0,\\
\left|a(x,z)-a(y,z)\right|^2&\leq& \ell_a|x-y|^2|z|^2.
\end{eqnarray*}

\begin{remark}

An example of the function pair $(a,b)$ is
$a(x,z)=b(x,z)=\sigma(x)z$ in $(H1)$, where $\sigma:\R\to\R$ is a
Lipschitzian map with Lip-coefficient $\sqrt{\ell_a}$ and
$\sigma(0)=0$. In the case, the perturbation in $(1.3)$ can be
rewritten as
$$\sigma(u(t))\D L_t,$$ where $(L_t)_{t\geq0}$ is a L\'{e}vy process
{\rm(}with L\'{e}vy measure $\pi(\cdot))$ given by
$$L_t=\int_0^t\!\!\int_{Z_1}z\widetilde{N}(\D z,\D s)+\int_0^t\!\!\int_{Z\setminus Z_1}z{N}(\D z,\D s),$$
by employing the L\'{e}vy-Khintchine Theorem $($see e.g. Sato
{\rm\cite{Sato}}$)$.
\end{remark}

In the coming section, we shall prove existence and uniqueness of
the weak solutions to (1.3). A $V\times H$-valued
$(\bar{\F}_t)_{t\geq0}$-adapted process
$X=(X(t))_{t\geq0}=((u(t),v(t)))_{t\geq0}$ is called a weak solution
of (1.3) with an initial value $X(0)=(\varphi,\psi)\in
V\times H$, if it fulfills the following two conditions:\\

\noindent{(1)}\quad $X\in C([0,T];V)\times
\DD([0,T];H)$\footnote{For $T>0$, $\DD([0,T];H)$ denotes the space
of all RCLL $(\bar{\F}_t)_{t\geq0}$-adapted random processes.} for
each $T>0$, $\Y$-a.s.\  \ and

\noindent{(2)}\quad For all test pairs $\phi=(\phi_1,\phi_2)^\T\in
D(\A^*)$, it holds that
\begin{eqnarray}
\left<X^\T(t),\phi\right>=\left<X^\T(0),\phi\right>+\int_0^t\left<
X^\T(s),\A^*\phi\right>\D s+\int_0^t\left<G(X^\T(s)),\phi\right>\D
s,
\end{eqnarray}
almost surely for $t\geq0$, where $X^\T(t)=(u(t),v(t))^\T$ and
$\A^*$ denotes the adjoint operator of $\A$ and $D(\A^*)$ is its
domain of the definition. In addition,
\begin{eqnarray*}
\A&=&\left[\begin{array}{cc} 0 & I \\ \J A\ & \J\kappa I
\end{array}\right],\nonumber\\
G(X^T(t))&=&\left[\begin{array}{c} 0 \\
\int_{Z_1}a(u(t\J),z)\dot{\widetilde{N}}(\D z,t)+\int_{Z\setminus
Z_1}b(u(t\J),z)\dot{N}(\D z,t)\end{array}\right].
\end{eqnarray*}

\section{Existence and uniqueness}
\setcounter{equation}{0}

The aim of this section is to establish the existence of a unique
weak solution for (1.3) under the condition $(H1)$.

The following result concentrates on the counterpart with small
jumps.
\begin{lemma}
Suppose that $h\in \LL^2([0,T]\times Z_1;V)$ and
$Y(0)=(\varphi,\psi)\in V\times H$. Then for any $T>0$, there exists
a unique weak solution $(Y(t))_{t\geq0}=((u(t),v(t)))_{t\geq0}\in
C([0,T];V)\times \DD([0,T];H)$ for the system:
\begin{eqnarray}
\left\{\begin{array}{ll} \D u(t)=v(t)\D t\\
\D v(t)=-\left[\kappa v(t)+A u(t)\right]\D
t+\int_{Z_1}h(t\J,z)\tilde{N}(\D z,\D t),\\
u(0)=\varphi,\ \ v(0)=\psi.
\end{array}\right.
\end{eqnarray}
\end{lemma}
\prf We are first to define,
\begin{eqnarray*}
g(t)=\int_0^t\!\!\int_{Z_1}h(s,z)\tilde{N}(\D z,\D s),\ \ \ \ \
t\geq0.
\end{eqnarray*}
Since $h\in \LL^2([0,T]\times Z_1;V)$, $g\in \LL^2([0,T];V)$. Let's
consider the system,
\begin{eqnarray}
\left\{\begin{array}{ll} \D u (t)=\left[\bar{v}(t)+g(t)\right]\D t\\
\D\bar{v}(t)=-[\kappa (\bar{v}(t)+g(t))+Au(t)]\D t
\end{array}\right.
\end{eqnarray}
In light of Lions \cite{JLL}, (3.2) admits a unique weak solution
$Z(t)=(u(t),\bar{v}(t))$ such that $Z\in C([0,T];V)\times
C([0,T];H)$. Let $v(t)=\bar{v}(t)+g(t)$. Then $Y(t)=(u(t),v(t))$
solves (3.1) and furthermore $Y\in C([0,T];V)\times \DD([0,T];H)$.
Thus we complete the proof of the lemma. \prfe

\begin{proposition}
Let the condition ${(H1})$ hold. Then for $X(0)=(\varphi,\psi)\in
V\times H$, there exists a unique weak solution
$X=(X(t))_{t\geq0}=((u(t),v(t)))_{t\geq0}$ for the system:
\begin{eqnarray}
\left\{\begin{array}{ll} \D u(t)=v(t)\D t\\
\D v(t)=-\left[\kappa v(t)+A u(t)\right]\D
t+\int_{Z_1}a(u(t\J),z)\tilde{N}(\D z,\D t),\\
u(0)=\varphi,\ \ v(0)=\psi.
\end{array}\right.
\end{eqnarray}
\end{proposition}

\prf Let's construct a sequence of $(\bar{\F}_t)_{t\geq0}$-adapted
random processes $(X^n)_{n\geq0}$ by $X^0(t)=X(0)$ for all $t\geq0$,
and for $n\geq0$,
$X^{n+1}=(X^{n+1}(t))_{t\geq0}=((u^{n+1}(t),v^{n+1}(t))_{t\geq0}\in
C([0,T];V)\times\DD([0,T];H)$ being the unique weak solution for the
following system:
\begin{eqnarray}
\left\{\begin{array}{ll} \D u^{n+1}(t)=v^{n+1}(t)\D t\\
\D v^{n+1}(t)=-\left[\kappa v^{n+1}(t)+A u^{n+1}(t)\right]\D
t+\int_{Z_1}a(u^n(t\J),z)\tilde{N}(\D z,\D t),\\
u^{n+1}(0)=\varphi,\ \ v^{n+1}(0)=\psi.
\end{array}\right.
\end{eqnarray}
By virtue of Lemma 3.1, it follows that $X^{n+1}$ exists. In what
follows, we show that the sequence $(X^n)_{n\geq1}$ is cauchy in
$C([0,T];V)\times \DD([0,T];H)$ compatibled with the uniform
topology. The It\^{o} rule (see e.g. Ikeda and Watanabe \cite{IW})
for $\left|v^{n+1}(t)-v^n(t)\right|^2$ yields that,
\begin{eqnarray}
&&\left|X^{n+1}(t)-X^n(t)\right|_{V\times H}^2\nonumber\\
&&\quad=\left\|u^{n+1}(t)-u^n(t)\right\|^2+\left|v^{n+1}(t)-v^n(t)\right|^2\nonumber\\
&&\quad=\left\|u^{n+1}(t)-u^n(t)\right\|^2-2\kappa\int_0^t\left|v^{n+1}(s)-v^n(s)\right|^2\D
s-\left\|u^{n+1}(t)-u^n(t)\right\|^2
\nonumber\\
&&\quad\quad+2\int_0^t\!\!\int_{Z_1}\left|a(u^{n}(s),z)-a(u^{n-1}(s),z)\right|^2
\pi(\D z)\D s\nonumber\\
&&\quad\quad+\int_0^t\!\!\int_{Z_1}[|(v^{n+1}(s\J)-v^n(s\J))
+(a(u^{n}(s\J),z)-a(u^{n-1}(s\J),z))|^2\nonumber\\
&&\quad\quad\quad\quad\quad\quad-|v^{n+1}(s\J)-v^n(s\J)|^2]\tilde{N}(\D
z,\D s)\nonumber\\
&&\quad=-2\kappa\int_0^t\left|v^{n+1}(s)-v^n(s)\right|^2\D
s+2\int_0^t\!\!\int_{Z_1}\left|a(u^{n}(s),z)-a(u^{n-1}(s),z)\right|^2
\pi(\D z)\D s\nonumber\\
&&\quad\quad+\int_0^t\!\!\int_{Z_1}[|(v^{n+1}(s\J)-v^n(s\J))
+(a(u^{n}(s\J),z)-a(u^{n-1}(s\J),z))|^2\nonumber\\
&&\quad\quad\quad\quad\quad\quad-|v^{n+1}(s\J)-v^n(s\J)|^2]\tilde{N}(\D
z,\D s).
\end{eqnarray}
In light of the condition $(H1)$ and Poincare-type inequality (2.5),
one gets,
\begin{eqnarray}
&&2\int_0^t\!\!\int_{Z_1}\left|a(u^{n}(s),z)-a(u^{n-1}(s),z)\right|^2
\pi(\D z)\D s\nonumber\\
&&\quad\leq
2\ell_a\int_0^t\!\!\int_{Z_1}|u^n(s)-u^{n-1}(s)|^2|z|^2\pi(\D
z)\D s\nonumber\\
&&\quad=\frac{2\bar{\theta}\ell_a}{\lambda_1}\int_0^t\left\|u^n(s)-u^{n-1}(s)\right\|^2\D
s.
\end{eqnarray}
Now we turn to the last term of the r.h.s. of (3.5). For $t\geq0$,
define
\begin{eqnarray}
II(t)&=&2\int_0^t\!\!\int_{Z_1}\left<v^{n+1}(s\J)-v^n(s\J),a(u^{n}(s\J),z)-a(u^{n-1}(s\J),z)\right>\tilde{N}(\D
z,\D s)\nonumber\\
&&\quad+\int_0^t\!\!\int_{Z_1}\left|a(u^{n}(s\J),z)-a(u^{n-1}(s\J),z)\right|^2
\tilde{N}(\D z,\D s)\nonumber\\
&:=&II_1(t)+II_2(t).
\end{eqnarray}
Then for the term $II_1$,
\begin{eqnarray}
&&[II_1,II_1]_t^{1/2}\nonumber\\
&&\quad=2\left[\int_0^t\!\!\int_{Z_1}\left<v^{n+1}(s\J)-v^n(s\J),a(u^n(s\J),z)-a(u^{n-1}(s\J),z)\right>^2N(\D
z,\D s)\right]^{1/2}\nonumber\\
&&\quad\leq2\left[\int_0^t\!\!\int_{Z_1}\left|v^{n+1}(s\J)-v^n(s\J)\right|^2
\left|a(u^n(s\J),z)-a(u^{n-1}(s\J),z)\right|^2N(\D
z,\D s)\right]^{1/2}\nonumber\\
&&\quad\leq2\sup_{0\leq s\leq
t}\left|v^{n+1}(s)-v^n(s)\right|\nonumber\\
&&\quad\quad\times\left[\int_0^t\!\!\int_{Z_1}\left|a(u^n(s\J),z)-a(u^{n-1}(s\J),z)\right|^2N(\D
z,\D s)\right]^{1/2}\nonumber\\
&&\quad\leq\frac{1}{4\sqrt{6}}\sup_{0\leq s\leq
t}\left|v^{n+1}(s)-v^n(s)\right|^2\nonumber\\
&&\quad\quad+4\sqrt{6}\int_0^t\!\!\int_{Z_1}\left|a(u^n(s\J),z)
-a(u^{n-1}(s\J),z)\right|^2N(\D z,\D s).
\end{eqnarray}
As a consequence, the Davis inequality and Poincare-type inequality
(2.5) jointly imply that,
\begin{eqnarray}
&&\E\left[\sup_{0\leq s\leq t}\left|II_1(s)\right|\right]\nonumber\\
&&\quad\leq2\sqrt{6}\E\left[[II_1,II_1]_t^{1/2}\right]\nonumber\\
&&\quad\leq\frac{1}{2}\E\left[\sup_{0\leq s\leq
t}\left|v^{n+1}(s)-v^n(s)\right|^2\right]
+\frac{48\bar{\theta}\ell_a}{\lambda_1}\int_0^t\E\left\|u^n(s)-u^{n-1}(s)\right\|^2\D
s.\nonumber\\
\end{eqnarray}
As for the term $II_2$, analogously we have,
\begin{eqnarray}
&&[II_2,II_2]_t^{1/2}\nonumber\\
&&\quad=\left[\int_0^t\!\!\int_{Z_1}\left|a(u^n(s\J),z)-a(u^{n-1}(s\J),z)\right|^4N(\D
z,\D s)\right]^{1/2}\nonumber\\
&&\quad
\leq\ell_a\left[\int_0^t\!\!\int_{Z_1}\left|u^n(s\J)-u^{n-1}(s\J)\right|^4z^4N(\D
z,\D s)\right]^{1/2}\nonumber\\
&&\quad\leq\frac{1}{16\sqrt{6}}\sup_{0\leq s\leq
t}\left\|u^n(s)-u^{n-1}(s)\right\|^2\nonumber\\
&&\quad\quad+\frac{4\sqrt{6}\ell_a^2}{\lambda_1^2}\int_0^t\!\!\int_{Z_1}
\left\|u^n(s\J)-u^{n-1}(s\J)\right\|^2z^4N(\D z,\D s),
\end{eqnarray}
and so,
\begin{eqnarray}
&&\E\left[\sup_{0\leq s\leq t}|II_2(s)|\right]\nonumber\\
&&\quad\leq\frac{1}{8}\E\left[\sup_{0\leq s\leq
t}\left\|u^{n}(s)-u^{n-1}(s)\right\|^2\right]
+\frac{48\bar{\theta}\ell_a^2}{\lambda_1^2}\int_0^t\E\left\|u^n(s)-u^{n-1}(s)\right\|^2\D
s,\nonumber\\
\end{eqnarray}
where we used the fact $\int_{Z_1}|z|^4\pi(\D z)\leq\bar{\theta}$.

In the following, we divide (3.5) into two respective parts
$\left\|u^{n+1}(t)-u^n(t)\right\|^2$ and
$\left|v^{n+1}(t)-v^n(t)\right|^2$ and estimate them respectively.
According to (3.5) and (3.6), we can conclude that for all $t>0$,
\begin{eqnarray*}
&&\E\left[\sup_{0\leq s\leq
t}\left\|u^{n+1}(s)-u^n(s)\right\|^2\right]+\E\left[\sup_{0\leq
s\leq
t}\left|v^{n+1}(s)-v^n(s)\right|^2\right]\nonumber\\
&&\quad\leq\frac{2\bar{\theta}\ell_a}{\lambda_1}\E\int_0^t\left\|u^{n}(s)-u^{n-1}(s)\right\|^2\D
s+\E\left[\sup_{0\leq s\leq t}II_1(s)\right]+\E\left[\sup_{0\leq
s\leq t}II_2(s)\right].
\end{eqnarray*}
From (3.9) and (3.11), it follows that,
\begin{eqnarray*}
&&\E\left[\sup_{0\leq s\leq
t}\left\|u^{n+1}(s)-u^n(s)\right\|^2\right]+\E\left[\sup_{0\leq
s\leq
t}\left|v^{n+1}(s)-v^n(s)\right|^2\right]\nonumber\\
&&\quad\leq\frac{1}{8}\E\left[\sup_{0\leq s\leq
t}\left\|u^{n}(s)-u^{n-1}(s)\right\|^2\right]+C_1\E\int_0^t\left\|u^{n}(s)-u^{n-1}(s)\right\|^2\D
s\nonumber\\
&&\quad\quad+\frac{1}{2}\E\left[\sup_{0\leq s\leq
t}\left|v^{n+1}(s)-v^n(s)\right|^2\right],
\end{eqnarray*}
where
$C_1=\frac{50\bar{\theta}\ell_a\lambda_1+48\bar{\theta}\ell_a^2}{\lambda_1^2}$.
This implies that
\begin{eqnarray*}
\E\left[\sup_{0\leq s\leq
t}\left\|u^{n+1}(s)-u^n(s)\right\|^2\right]&\leq&\frac{1}{8}\E\left[\sup_{0\leq
s\leq t}\left\|u^n(s)-u^{n-1}(s)\right\|^2\right]\nonumber\\
&&+C_1 \E\int_0^t\left\|u^n(s)-u^{n-1}(s)\right\|^2\D s.
\end{eqnarray*}
Analogously, using (3.5) and (3.6), one gets,
\begin{eqnarray*}
&&\E\left[\sup_{0\leq s\leq
t}\left\|v^{n+1}(s)-v^n(s)\right\|^2\right]\nonumber\\
&&\quad\leq-2\kappa\E\int_0^t\left|v^{n+1}(s)-v^n(s)\right|^2\D
s+\frac{2\bar{\theta}\ell_a}{\lambda_1}\E\int_0^t\left\|u^{n}(s)-u^{n-1}(s)\right\|^2\D
s\nonumber\\
&&\quad\quad+\E\left[\sup_{0\leq s\leq
t}II_1(s)\right]+\E\left[\sup_{0\leq s\leq t}II_2(s)\right].
\end{eqnarray*}
We also apply (3.9) and (3.11) to conclude that
\begin{eqnarray*}
&&\E\left[\sup_{0\leq s\leq
t}\left|v^{n+1}(s)-v^n(s)\right|^2\right]\nonumber\\
&&\quad\leq-2\kappa\E\int_0^t\left|v^{n+1}(s)-v^n(s)\right|^2\D
s+C_1\E\int_0^t\left\|u^{n}(s)-u^{n-1}(s)\right\|^2\D
s\nonumber\\
&&\quad\quad+\frac{1}{8}\E\left[\sup_{0\leq s\leq
t}\left\|u^n(s)-u^{n-1}(s)\right\|^2\right]+\frac{1}{2}\E\left[\sup_{0\leq
s\leq t}\left|v^{n+1}(s)-v^n(s)\right|^2\right].
\end{eqnarray*}
As a consequence,
\begin{eqnarray*}
&&\E\left[\sup_{0\leq s\leq
t}\left|v^{n+1}(s)-v^n(s)\right|^2\right]\nonumber\\
&&\quad\leq-4\kappa\E\int_0^t\left|v^{n+1}(s)-v^n(s)\right|^2\D
s+\frac{1}{4}\E\left[\sup_{0\leq s\leq
t}\left\|u^n(s)-u^{n-1}(s)\right\|^2\right]\nonumber\\
&&\quad\quad+2C_1\E\int_0^t\left\|u^n(s)-u^{n-1}(s)\right\|^2\D s.
\end{eqnarray*}
Consequently, for all $t>0$,
\begin{eqnarray}
&&\E\left[\sup_{0\leq s\leq
t}\left|X^{n+1}(s)-X^n(s)\right|_{V\times H}^2\right]\nonumber\\
&&\quad\leq-4\kappa\E\int_0^t\left|v^{n+1}(s)-v^n(s)\right|^2\D
s+\frac{3}{8}\E\left[\sup_{0\leq s\leq
t}\left|X^n(s)-X^{n-1}(s)\right|_{V\times H}^2\right]\nonumber\\
&&\quad\quad+3C_1\E\int_0^t\left|X^n(s)-X^{n-1}(s)\right|_{V\times
H}^2\D s.
\end{eqnarray}
For each $0<t\leq T$, let $V^n(t)=\E\left[\sup_{0\leq s\leq
t}\left|X^{n+1}(s)-X^n(s)\right|_{V\times H}^2\right]$ with
$n\geq0$. Then (3.12) can be rewritten as
\begin{eqnarray*}
V^n(t)&\leq&\frac{3}{8}V^{n-1}(t)+3C_1\int_0^tV^{n-1}(s)\D s,\ \ \ \
n\geq1,
\end{eqnarray*}
A recursive scheme for the above relation between $V^n$ and
$V^{n-1}$ shows that for each $T>0$, there exists a constant $C_T>0$
such that
\begin{eqnarray*}
V^n(t)&\leq& C_T\sum_{i=0}^n{\rm
C}_n^{i}(\frac{3}{8})^{n-i}\frac{C_T^i}{i!}=C_T(\frac{3}{8})^{n}\sum_{i=0}^n{\rm
C}_n^{i}\frac{(8C_T/3)^i}{i!}\nonumber\\
&\leq&C_T(\frac{3}{4})^n\exp\left(\frac{8C_T}{3}\right),
\end{eqnarray*}
where we used the fact $\sum_{i=0}^n{\rm C}_n^{i}=2^n$ and hence
${\rm C}_n^{i}\leq 2^n$ for each $i=0,1,\dots,n$. This recursive
result further yields that there exists a random process $X\in
C([0,T];V)\times\DD([0,T];H)$ such that
\begin{eqnarray}
\lim_{n\to\infty}\E\left[\sup_{0\leq t\leq
T}\left|X^n(t)-X(t)\right|_{V\times H}^2\right]=0.
\end{eqnarray}
Letting $n\to+\infty$ in (3.4) to conclude that  $(X(t))_{t\geq0}$
is a weak solution of (3.3). The uniqueness of $(X(t))_{t\geq0}$
follows from the It\^{o} rule and Gronwall Lemma. We omit its proof.
\prfe

\begin{theorem}
Suppose that the condition ${(H1)}$ holds. Then for
$X(0)=(\varphi,\psi)\in V\times H$, $(1.3)$ admits a unique weak
solution $X=(X(t))_{t\geq0}=(u(t),v(t))_{t\geq0}$.
\end{theorem}

\prf It follows from (1.1) that, $\pi(Z\setminus Z_1)<\infty$. Then
the process $(N(Z\setminus Z_1\times[0,t]))_{t\geq0}$ has only
finite jumps in each finite interval of $\R_+$, i.e., there exist
increasing jump times $0<\tau_1<\tau_2<\cdots<\tau_n<\cdots$.
Moreover, $(N(A\times[0,t]))_{(A,t)\in\B(Z\setminus Z_1)\times\R_+}$
can be represented by a $Z$-valued point process $(\q(t))_{t\geq0}$
with the domain $D_{\q}$ as a countable subset of $\R_+$. That is,
\begin{eqnarray}
N(A\times[0,t])=\sum_{s\in D_\q,s\leq t}{\bf1}_A(\q(s)),\ \ {\rm
for}\ t>0\ {\rm and}\ A\in{\mathcal{B}}(Z\setminus Z_1).
\end{eqnarray}
Therefore for $k=1,2,\dots$, $\tau_k\in\{t\in D_{\q};\ \q(t)\in
Z\setminus Z_1\}$. For each $n\in\N$, we easily see that $\tau_k$ is
an $(\bar{\F}_t)_{t\geq0}$-stopping time and $\tau_k\to\infty$, as
$k\to\infty$. For each $T\in(0,\tau_1)$, By virtue of Proposition
3.2, there exists a unique weak solution $X^0\in C([0,T];V)\times
\DD([0,T];H)$ on $[0,\tau_1)$. Construct the following
\begin{eqnarray*}
X^1(t)=\left \{\begin{array}{ll}
           X^0(t), &\ \ t\in[0,\tau_1),\\
           X^0({\tau_1\J})+
           \left[\begin{array}{c} 0 \\
                  b\left(u({\tau_1\J}),\q(\tau_1)\right)
                  \end{array}\right]^\T, & \quad t=\tau_1.
        \end{array}
\right.
\end{eqnarray*}
Therefore $(X^1(t))_{0\leq t\leq\tau_1}$ uniquely solves (3.1) in
the time interval $[0,\tau_1]$. Furthermore we define
\begin{eqnarray*}
\left \{\begin{array}{cl}
           \tilde{X}_0^1 & =X^1({\tau_1}),\\
           \widetilde{\q}(t) & =\q(t+\tau_1),\\
           D_{\tilde{\q}} & =\left\{t\geq0;\ t+\tau_1\in
           D_\q\right\},\\
           \widetilde{{\F}}_t & =\bar{\F}_{\tau_1+t}.
        \end{array}
\right.
\end{eqnarray*}
Note that $\tau_2-\tau_1\in \{t\in D_{\tilde{\q}};\
\widetilde{\q}(t)\in Z\setminus Z_1\}$. Then we can construct a
process $(\widetilde{X}^1(t))_{0\leq t\leq\tau_2-\tau_1}$ by a same
way as for $(X^1(t))_{0\leq t\leq\tau_1}$. Thus we let
\begin{eqnarray*}
X^2(t)=\left \{\begin{array}{ll}
           X^1(t), &\quad 0\leq t\leq\tau_1,\\
           \widetilde{X}^1({t-\tau_1}), &\quad \tau_1\leq
           t\leq\tau_2.
        \end{array}
\right.
\end{eqnarray*}
Then $X^2(t)$ is a unique weak solution of (1.3) in the time
interval $[0,\tau_2]$. Hence the existence of the unique global weak
solution follows from the above successive procedure, and the
theorem is proved. \prfe

\section{Invariant measure}
\setcounter{equation}{0}

In the section, we shall study the existence of a unique invariant
measure associated with the transient semigroup $(\p_t)_{t\geq0}$
defined by
\begin{eqnarray}
\p_t\Phi((\varphi,\psi))=\E\left[\Phi(X_t^0((\varphi,\psi)))\right],\
\ (\varphi,\psi)\in V\times H,\ \ \Phi\in C_b(V\times H),
\end{eqnarray}
where $X_t^0((\varphi,\psi))=(u_t^0(\varphi),v_t^0(\psi))$ denotes
the weak solution of (1.3) with the initial value $(\varphi,\psi)\in
V\times H$ at time-zero. As for the Markov property of
$X_t^0((\varphi,\psi))$, we refer to Bo et al. \cite{BSW08}.

To establish the invariant measure for $(\p_t)_{t\geq0}$, set
\begin{eqnarray}
\delta_0=\frac{\lambda_1}{2\kappa}\wedge\frac{\kappa}{4},
\end{eqnarray}
and $\rho_\delta(t)=\delta u(t)+v(t)$ with the weak solution
$(X(t))_{t\geq0}=(u(t),v(t))_{t\geq0}$ to (1.3). Then we claim that,
\begin{lemma}
For all positive $\delta\leq\delta_0$ and $t\geq0$, it holds that
\begin{eqnarray}
\left|\rho_\delta(t)\right|^2+\left\|u(t)\right\|^2&\leq&\left|\delta\varphi+\psi\right|^2+\left\|\varphi\right\|^2
-\int_0^t\left[\delta\left\|u(s)\right\|^2
+\kappa\left|\rho_\delta(s)\right|^2\right]\D s\nonumber\\
&&+\int_0^t\!\!\int_{Z_1}\left|a(u(s),z)\right|^2\pi(\D z)\D s+M_t\nonumber\\
&&+\int_0^t\!\!\int_{Z\setminus
Z_1}\left[\left|b(u(s),z)\right|^2+2\left<\rho_\delta(t),b(u(s),z)\right>\right]\pi(\D
z)\D s,\nonumber\\
\end{eqnarray}
where $(M_t)_{t\geq0}$ is a {\rm RCLL}
$(\bar{\F}_t)_{t\geq0}$-martingale with mean zero and which is given
by
\begin{eqnarray*}
M_t&=&\int_0^t\!\!\int_{Z_1}\left[\left|\rho_\delta(s\J)+a(u(s\J),z)\right|^2-\left|\rho_\delta(s\J)\right|^2\right]
\widetilde{N}(\D z,\D s)\nonumber\\
&&+\int_0^t\!\!\int_{Z\setminus
Z_1}\left[\left|\rho_\delta(s\J)+b(u(s\J),z)\right|^2-\left|\rho_\delta(s\J)\right|^2\right]
\widetilde{N}(\D z,\D s),\ \ t\geq0.
\end{eqnarray*}
\end{lemma}

\prf By virtue of (1.3), the process $(\rho_\delta(t))_{t\geq0}$ is
a RCLL $(\bar{\F}_t)_{t\geq0}$-semimartingale which satisfies the
following dynamics,
\begin{eqnarray}
\D\rho_\delta(t)&=&(\delta-\kappa)\rho_\delta(t)\D
t-\left[\delta(\delta-\kappa)+A\right]u(t)\D
t+\int_{Z_1}a(u(t\J),z)\widetilde{N}(\D z,\D t)\nonumber\\
&&+\int_{Z\setminus Z_1}b(u(t\J),z)N(\D z,\D t),\\
\rho_\delta(0)&=&\delta\varphi+\psi.\nonumber
\end{eqnarray}
On the other hand, we remark that for $\delta\leq\delta_0$ and
$t\geq0$,
\begin{eqnarray}
&&\delta(\kappa-\delta)\left<u(t),\rho_\delta(t)\right>-(\kappa-\delta)\left|\rho_\delta(t)\right|^2
-\delta\left\|u(t)\right\|^2\nonumber\\
&&\quad\quad\leq-\frac{\delta}{2}\left\|u(t)\right\|^2-\frac{\kappa}{2}\left|\rho_\delta(t)\right|^2.
\end{eqnarray}
Then apply the It\^{o} rule w.r.t. Poisson random measures (see
Ikeda and Watanabe \cite{IW}) to
$\frac{1}{2}\left|\rho_\delta(t)\right|^2$, the desired result
follows from (4.4) and (4.5) immediately. \prfe

Hereafter, we define an energy functional $\EE^\delta$ on $V\times
H$ by
\begin{eqnarray*}
\EE^\delta(u,v)=|\delta u+v|^2+\|u\|^2,\ \ \ (u,v)\in V\times H.
\end{eqnarray*}
In order to explore the invariant measure, we impose the following
condition on the function $b:\R\times Z\setminus Z_1\to\R$,

\noindent$({H2})$\quad There exists $\ell_b>0$ such that
\begin{eqnarray*}
b(0,z)&\equiv& 0,\\
|b(x,z)-b(y,z)|^2&\leq&\ell_b|x-y|^2.
\end{eqnarray*}
\begin{remark}
Note that the condition $(H2)$ rules out the case of
$b(x,z)=\sigma(x)z$ in Remark $2.1$. To incorporate the case into
the section, we impose the condition,

\noindent$(H2)'$\quad There exists $\ell_b>0$ such that
\begin{eqnarray*}
b(0,z)&\equiv& 0,\\
|b(x,z)-b(y,z)|^2&\leq&\ell_b|x-y|^2|z|^{p},\ \ {\rm with\ the\
integer\ }p\geq2,\\
\theta_p&=&\int_{Z\setminus Z_1}|z|^p\pi(\D z)<\infty.
\end{eqnarray*}
The last condition in $(H2)'$ is equivalent to that the L\'{e}vy
process $(L_t)_{t\geq0}$ admits the finite $p$-order moment.
Compared with $(H2)$ and $(H2)'$, we also note that if $(H2)$ holds,
then L\'{e}vy measure $\pi(\cdot)$ is unrestrictive. However it
rules out the case in Remark $2.1$. If $(H2)'$ is assumed to be
true, then the case in Remark $2.1$ is included, but an additional
condition on $\pi(\cdot):\ \theta_2<\infty$ has to be imposed.
However the essential proofs in the section by employing $(H2)$ and
$(H2)'$ are indistinctive.
\end{remark}
Consequently,
\begin{lemma}
Suppose the triple $(\ell_a,\ell_b,\kappa)$ satisfies that,
\begin{eqnarray}
\frac{\bar{\theta}\ell_a+2\underline{\theta}\ell_b}
{\lambda_1}<\delta_0,\ \ {\rm and}\ \ \kappa>{\underline{\theta}},
\end{eqnarray}
where $\bar{\theta}$, $\underline{\theta}$ are defined in $(1.2)$.
Then under the conditions $(H1)$--$(H2)$, or under the conditions
$(H1)$--$(H2)'$ for the triple $(\ell_a,\ell_b,\kappa)$ satisfying
$(4.6)$ with $\underline{\theta}$ replaced by $\theta_p$, there
exist positive constants $\delta\leq\delta_0$ and
$\lambda=\lambda(\delta)$ such that
\begin{eqnarray*}
\EE^\delta(u(t),v(t))&\leq&\EE^\delta(\varphi,\psi)-\lambda\int_0^t\EE^\delta(u(s),v(s))\D
s+M_t,\ \ \ t\geq0,
\end{eqnarray*}
where the RCLL $(\bar{\F}_t)_{t\geq0}$-martingale $(M_t)_{t\geq0}$
is defined in Lemma $4.1$.
\end{lemma}

\begin{remark}
$1.$ Note that the parameter $\delta_0$ depends on $\kappa$ $($see
$(4.2))$. However, we can choose a pair
$(\ell_a^*,\ell_b^*)\in(0,\infty)^2$ $($at least when they are small
enough$)$ such that
\begin{eqnarray*}
\underline{\theta}\vee\sqrt{2\lambda_1}<
\frac{\lambda_1^2}{2\bar{\theta}\ell_a^*+4\underline{\theta}\ell_b^*}.
\end{eqnarray*}
Taking any $\kappa^*\in(\underline{\theta}\vee\sqrt{2\lambda_1},
{\lambda_1^2}/{[2\bar{\theta}\ell_a^*+4\underline{\theta}\ell_b^*]})$.
Then the triple $(\ell_a^*,\ell_b^*,\kappa^*)$ fulfills $(4.6)$.

\noindent$2.$ If the condition $(H2)$ is placed by $(H2)'$, then the
constant $\underline{\theta}$ should be placed by $\theta_p$ in
$(4.6)$. In the case, we choose a pair
$(\ell_a^*,\ell_b^*)\in(0,\infty)^2$ $($at least when they are small
enough$)$ such that
\begin{eqnarray*}
{\theta_p}\vee\sqrt{2\lambda_1}<
\frac{\lambda_1^2}{2\bar{\theta}\ell_a^*+4{\theta_p}\ell_b^*}.
\end{eqnarray*}
\end{remark}

We are now in a position to prove Lemma 4.3.

\noindent{\bf Proof of Lemma 4.3.}\quad Using the conditions
$(H1)$--$(H2)$ and Poincare-type inequality (2.5), it follows that
\begin{eqnarray}
\int_{Z_1}\left|a(u(t),z)\right|^2\pi(\D
z)\leq\bar{\theta}\ell_a\left|u(t)\right|^2\leq\frac{\bar{\theta}\ell_a}{\lambda_1}\left\|u(t)\right\|^2,\
\ t\geq0,
\end{eqnarray}
and
\begin{eqnarray}
&&\left|\int_{Z\setminus
Z_1}\left[\left|b(u(t),z)\right|^2+2\left<\rho_\delta(t),b(u(t),z)\right>\right]\pi(\D
z)\right|\nonumber\\
&&\quad \leq
\frac{2\underline{\theta}\ell_b}{\lambda_1}\|u(t)\|^2+{\underline{\theta}}\left|\rho_\delta(t)\right|^2,\
\ \ t\geq0.
\end{eqnarray}
Thanks to (4.6), we can choose a positive $\delta\in
({\bar{\theta}\ell_a}/{\lambda_1}+{2\underline{\theta}\ell_b}/{\lambda_1},\delta_0]$,
and then Lemma 4.1 yields that,
\begin{eqnarray*}
\left|\rho_\delta(t)\right|^2+\left\|u(t)\right\|^2&\leq&\left|\rho_\delta(0)\right|^2+\left\|\varphi\right\|^2
-\int_0^t[\delta-{\bar{\theta}\ell_a}/{\lambda_1}-{2\underline{\theta}\ell_b}/{\lambda_1}]
\left\|u(s)\right\|^2\D s\nonumber\\
&&-\int_0^t[\kappa-{\underline{\theta}}]\left|\rho_\delta(s)\right|^2\D
s+M_t\nonumber\\
&\leq&\left|\rho_\delta(0)\right|^2+\left\|\varphi\right\|^2-\lambda\int_0^t
[\left|\rho_\delta(s)\right|^2+\left\|u(s)\right\|^2]\D s+M_t,
\end{eqnarray*}
where $\lambda=\min\{\delta-{\bar{\theta}\ell_a}/{\lambda_1}
-{2\underline{\theta}\ell_b}/{\lambda_1},\kappa-{\underline{\theta}}\}>0$.
When the conditions $(H1)$--$(H2)'$ are satisfied, the estimates
(4.7) and (4.8) also hold with $\underline{\theta}$ replaced by
$\theta_p$. Thus the proof of the lemma is complete. \prfe

In what follows, we state the main result of the section.
\begin{theorem}
Under the same conditions as in Lemma $4.3$, there exists a unique
invariant measure $\nu(\cdot)$ on $(V\times H,\B(V\times H))$ for
the transient semigroup $(\p_t)_{t\geq0}$ defined by $(4.1)$.
\end{theorem}

\prf We adopt the method used in Chow \cite{Chow06}. Let
$(\bar{N}(A\times[0,t]))_{A\in\B(Z)}$ be an independent copy of the
Poisson random measure $({N}(A\times[0,t]))_{A\in\B(Z)}$ for
$t\geq0$. For any $A\in\B(Z)$ and $t\in\R$,  define
\begin{eqnarray*}
\hat{N}(A\times[0,t])&=&N(A\times[0,t]),\quad {\rm if}\ t\geq0,\
{\rm and}\\
\hat{N}(A\times[t,0])&=&\bar{N}(A\times[0,-t]), \quad {\rm if}\ t<0.
\end{eqnarray*}
Let $\widetilde{\hat{N}}$ be the compensated Poisson random measure
of $\hat{N}$. For each $s\in\R$, consider the system:
\begin{eqnarray}
\left\{\begin{array}{ll} \D u(t,\xi) = v(t,\xi)\D t,\\
\D v(t,\xi)= -[\kappa v(t,\xi)+A u(t,\xi)]\D t
+\int_{Z_1}a(u(t\J,\xi),z)\widetilde{\hat{N}}(\D z,\D
t)\\
\ \ \ \ \ \ \ \ \ \ \ \ +\int_{Z\setminus
Z_1}b(u(t\J,\xi),z){\hat{N}}(\D z,\D
t),\\
u(s,\xi)=\varphi(\xi),\ \ v(s,\xi)=\psi(\xi),\ \ \xi\in D.
\end{array}\right.
\end{eqnarray}
By virtue of Theorem 3.3, there exists a unique solution
$(X_t^s((\varphi,\psi)))_{t>s}\in C([s,T];V)\times\DD([s,T];H)$ for
each $T>0$, provided $(\varphi,\psi)\in V\times H$. Therefore, from
the Gronwall Lemma, it follows that for some positive constants
$\delta\leq\delta_0$ and $\lambda=\lambda(\delta)$,
\begin{eqnarray}
\E\left[\EE^\delta(X_t^s((\varphi,\psi)))\right]&\leq&e^{-\lambda(t-s)}
\E\left[\EE^\delta(\varphi,\psi)\right],\ \ \ \ t>s.
\end{eqnarray}
For $s_1>s_2>0$, define
\begin{eqnarray*}
\hat{X}_t^{1,2}((\varphi,\psi))=(\hat{u}(t),\hat{v}(t))
=\left(u_t^{-s_1}(\varphi)-u_t^{-s_2}(\varphi),v_t^{-s_1}(\psi)-v_t^{-s_2}(\psi)\right).
\end{eqnarray*}
Then $\hat{X}_t^{1,2}((\varphi,\psi))$ fulfills that
\begin{eqnarray}
\left\{\begin{array}{ll} \D \hat{u}(t,\xi) = \hat{v}(t,\xi)\D t,\\
\D\hat{v}(t,\xi)= -[\kappa \hat{v}(t,\xi)+A\hat{u}(t,\xi)]\D t
+\int_{Z_1}\hat{a}(u_{t\J}^{\J s_1}(\xi),u_{t\J}^{\J
s_2}(\xi),z)\widetilde{\hat{N}}(\D z,\D
t)\\
\ \ \ \ \ \ \ \ \ \ \ \ +\int_{Z\setminus Z_1}\hat{b}(u_{t\J}^{\J
s_1}(\xi),u_{t\J}^{\J s_2}(\xi),z){\hat{N}}(\D z,\D
t),\\
\hat{u}(\J s_2,\xi)=u_{\J s_2}^{\J s_1}(\xi)-\varphi(\xi),\ \
\hat{v}(\J s_1,\xi)=v_{\J s_2}^{\J s_1}(\xi)-\psi(\xi),
\end{array}\right.
\end{eqnarray}
where for $(\xi,z)\in D\times Z$,
\begin{eqnarray*}
\hat{a}(u_{t}^{\J s_1}(\xi),u_{t}^{\J s_2}(\xi),z)&:=&a(u_{t}^{\J
s_1}(\xi),z)-a(u_{t}^{\J s_2}(\xi),z),\\
\hat{b}(u_{t}^{\J s_1}(\xi),u_{t}^{\J s_2}(\xi),z)&:=&b(u_{t}^{\J
s_1}(\xi),z)-b(u_{t}^{\J s_2}(\xi),z).
\end{eqnarray*}
Let $\hat{\rho}(t)=\delta\hat{u}(t)+\hat{v}(t)$ with $t\geq0$. Then
from Lemma 4.3, it follows that there exist positive
$\delta\leq\delta_0$ and $\lambda=\lambda(\delta)$ such that
\begin{eqnarray}
\E\left[\EE^\delta(\hat{X}_t^{1,2}((\varphi,\psi)))\right]&\leq&e^{-\lambda(t+s_2)}
\E\left[\EE^\delta(\hat{u}(\J s_2),\hat{v}(\J s_2))\right],\ \ \
t>\J s_2.
\end{eqnarray}
Thanks to (4.10), there exists a positive constant $C>0$ such that
\begin{eqnarray}
\E\left[\EE^\delta(\hat{X}_t^{1,2}((\varphi,\psi)))\right]&\leq&
Ce^{-\lambda(t+s_2)}[1+\EE^\delta(\varphi,\psi)] ,\ \ \ t>\J s_2.
\end{eqnarray}
Then by virtue of (4.13), one gets,
\begin{eqnarray}
\E\left[\EE^\delta(X_{0}^{\J s_1}((\varphi,\psi))-X_{0}^{\J
s_2}((\varphi,\psi)))\right]&\leq& Ce^{-\lambda
s_2}[1+\EE^\delta(\varphi,\psi)].
\end{eqnarray}
This implies that $(X_{0}^{\J s})_{s\geq0}$ is Cauchy in
$L^2(\Omega;V\times H)$. As a consequence, there exists a unique
random vector $X_0^{-\infty}((\varphi,\psi))\in L^2(\Omega;V\times
H)$ such that $X_0^{-s}((\varphi,\psi))\to
X_0^{-\infty}((\varphi,\psi))$, as $s\to\infty$ in
$L^2(\Omega;V\times H)$ sense. We remark that the vector processes
\begin{eqnarray*}
X_{0}^{\J s}((\varphi,\psi))=(u_0^{\J s}(\varphi),v_0^{\J s}(\psi))\
\ {\rm and}\ \
X_{s}^{0}((\varphi,\psi))=(u_s^{0}(\varphi),v_s^{0}(\psi))
\end{eqnarray*}
admit the same distribution on the same probability space for each
$s\geq0$. Let $\nu(\cdot)$ be the induced probability measure of
$X_0^{-\infty}((\varphi,\psi))$ on $(V\times H,\B(V\times H))$. Then
$\nu(\cdot)$ is the unique invariant measure for the transient
semigroup $(\p_t)_{t\geq0}$. Thus the proof of the theorem is
finished. \prfe\\

\noindent{\bf Acknowledgements.} The authors would like to thank an
anonymous referee and an Associated Editor of the journal for their
valuable comments and suggestions.


\small
\begin{thebibliography}{}

\bibitem{BDT} Barbu, V., Da Prato, G. and Tubaro, L. (2007). Stochastic wave
equations with dissipative damping. {\it Stoch. Process. Appl.} {\bf
117}, 1001-1013.

\bibitem{BSW} Bo, L., Shi, K. and Wang, Y. (2007). On a nonlocal
stochastic Kuramoto-Sivashinsky equation with jumps. {\it Stoch.
Dyn.} {\bf7}, 439-457.

\bibitem{BSW08} Bo, L., Shi, K. and Wang, Y. (2008). Stochastic wave equation
driven by compensated Poisson random measure. {Preprint}.

\bibitem{BTW} Bo, L., Tang, D. and Wang, Y. (2008). Explosive solutions of stochastic wave
equations with damping on $\R^d$. {\it J. Differential Equations.}
{\bf 244}, 170-187.

\bibitem{BMS}  Brzeniak, Z., Maslowski, B. and Seidler, J. (2005). Stochastic nonlinear beam
equations. {\it Probab. Th. Relat. Fields} {\bf132}, 119-149.


\bibitem{Chow02} Chow, P. (2002). Stochastic wave equation with polynimial
nonlinearity. {\it Ann. Appl. Probab.} {\bf12}, 361-381.

\bibitem{Chow06} Chow, P. (2006). Asymptotics of solutions to semilinear
stochastic wave equations. {\it Ann. Appl. Probab.} {\bf16},
757-789.

\bibitem{DaZa} Da Prato, G. and Zabczyk, J. (1996). {\it Ergodicity of Infinite
Dimensional Systems}. in: {London Mathematical Society Lecture
Notes}, vol 229, Cambridge University Press.

\bibitem{IW} Ikeda, N. and Watanabe, S. (1981). {\it Stochastic Differential Equations and Diffusion
Processes}. North-Holland, Amsterdam.

\bibitem{JLL} Lions, J. L. (1961). {\it Equations differentielles operationelles et problemes aux
limits}. Springer-Verlag, Berlin.

\bibitem{PeZabc06} Peszat, S. and Zabczyk, J. (2006). Stochastic heat and wave
equations driven by an impulsive noise, in: Da Prato, G. and Tubaro,
L. (Eds.), {\it Stochastic Partial Differential Equations and
Applications--VII}, in: Lect. Notes Pure Appl., vol. 245, Chapman \&
Hall/CRC, Boca Raton, pp. 229-242.

\bibitem{PeZabc} Peszat, S. and  Zabczyk, J. (2007). {\it Stochastic Partial
Differential Equations with L\'{e}vy noise: An Evolution Equation
Approach}. in: {Encyclopedia of Mathematics and Its Applications},
vol 113, Cambridge University Press.

\bibitem{Robin} Robinson, J. (2001). {\it Infinite-Dimensional Dynamical Systems}.
Cambridge University Press.

\bibitem{Sato} Sato, K. (1999). {\it L\'{e}vy Processes and Infinitely
Divisible Distributions}. Cambridge University Press.

\bibitem{Sch} Schiff, L. (1951). Nonlinear meson theory of nuclear forces I.
{\it Phys. Rev.} {\bf84}, 1-9.

\bibitem{Segal} Segal, I. (1963). The global Cauchy problem for a relativistic scalar
field with power interaction. {\it Bull. Soc. Math. France} {\bf91},
129-135.

\bibitem{Temam} Temam, R. (1997). {\it Infinite-Dimensional Dynamical Systems in
Mechanics and Physics, second ed.} Springer-Verlag, New York.

\bibitem{Wa} Walsh, J. (1986). {\it An Introduction to Stochastic Partial Differential
Equations}. Lecture Notes in Math. vol. 1180, Springer, Berlin, pp.
265-439.

\bibitem{Zei} Zeidler, E. (1990). {\it Nonlinear Functional Analysis and Its Applications, II/B,
Nonlinear Monotone Operators}. Springer-Verlag, New York.






\end{thebibliography}
\end{document}